\documentclass[twoside,leqno]{article}
\usepackage{amsmath}
\usepackage{amsthm}
\usepackage{amssymb}
\usepackage{amscd}
\usepackage{graphicx}
\usepackage{epsfig}

\usepackage{latexsym}
\usepackage{amsfonts}

\input xy
\xyoption{all} \tolerance=500

=1in \textwidth 15.6cm \topmargin 15mm\textheight23cm \advance\topmargin by -\headheight\advance\topmargin by
-\headsep

\numberwithin{equation}{section} \allowdisplaybreaks

\newtheorem{proposition}{Proposition}[section]

\theoremstyle{definition}
\newtheorem{definition}{Definition}[section]

\begin{document}
\font\black=cmbx10 \font\sblack=cmbx7 \font\ssblack=cmbx5 \font\blackital=cmmib10  \skewchar\blackital='177
\font\sblackital=cmmib7 \skewchar\sblackital='177 \font\ssblackital=cmmib5 \skewchar\ssblackital='177
\font\sanss=cmss10 \font\ssanss=cmss8 
\font\sssanss=cmss8 scaled 600 \font\blackboard=msbm10 \font\sblackboard=msbm7 \font\ssblackboard=msbm5
\font\caligr=eusm10 \font\scaligr=eusm7 \font\sscaligr=eusm5 \font\blackcal=eusb10 \font\fraktur=eufm10
\font\sfraktur=eufm7 \font\ssfraktur=eufm5 \font\blackfrak=eufb10

\font\bsymb=cmsy10 scaled\magstep2
\def\all#1{\setbox0=\hbox{\lower1.5pt\hbox{\bsymb
       \char"38}}\setbox1=\hbox{$_{#1}$} \box0\lower2pt\box1\;}
\def\exi#1{\setbox0=\hbox{\lower1.5pt\hbox{\bsymb \char"39}}
       \setbox1=\hbox{$_{#1}$} \box0\lower2pt\box1\;}

\def\mi#1{{\fam1\relax#1}}
\def\tx#1{{\fam0\relax#1}}

\newfam\bifam
\textfont\bifam=\blackital \scriptfont\bifam=\sblackital \scriptscriptfont\bifam=\ssblackital
\def\bi#1{{\fam\bifam\relax#1}}

\newfam\blfam
\textfont\blfam=\black \scriptfont\blfam=\sblack \scriptscriptfont\blfam=\ssblack
\def\rbl#1{{\fam\blfam\relax#1}}

\newfam\bbfam
\textfont\bbfam=\blackboard \scriptfont\bbfam=\sblackboard \scriptscriptfont\bbfam=\ssblackboard
\def\bb#1{{\fam\bbfam\relax#1}}

\newfam\ssfam
\textfont\ssfam=\sanss \scriptfont\ssfam=\ssanss \scriptscriptfont\ssfam=\sssanss
\def\sss#1{{\fam\ssfam\relax#1}}

\newfam\clfam
\textfont\clfam=\caligr \scriptfont\clfam=\scaligr \scriptscriptfont\clfam=\sscaligr
\def\cl#1{{\fam\clfam\relax#1}}

\newfam\frfam
\textfont\frfam=\fraktur \scriptfont\frfam=\sfraktur \scriptscriptfont\frfam=\ssfraktur
\def\fr#1{{\fam\frfam\relax#1}}

\def\cb#1{\hbox{$\fam\gpfam\relax#1\textfont\gpfam=\blackcal$}}

\def\hpb#1{\setbox0=\hbox{${#1}$}
    \copy0 \kern-\wd0 \kern.2pt \box0}
\def\vpb#1{\setbox0=\hbox{${#1}$}
    \copy0 \kern-\wd0 \raise.08pt \box0}

\def\pmb#1{\setbox0\hbox{${#1}$} \copy0 \kern-\wd0 \kern.2pt \box0}
\def\pmbb#1{\setbox0\hbox{${#1}$} \copy0 \kern-\wd0
      \kern.2pt \copy0 \kern-\wd0 \kern.2pt \box0}
\def\pmbbb#1{\setbox0\hbox{${#1}$} \copy0 \kern-\wd0
      \kern.2pt \copy0 \kern-\wd0 \kern.2pt
    \copy0 \kern-\wd0 \kern.2pt \box0}
\def\pmxb#1{\setbox0\hbox{${#1}$} \copy0 \kern-\wd0
      \kern.2pt \copy0 \kern-\wd0 \kern.2pt
      \copy0 \kern-\wd0 \kern.2pt \copy0 \kern-\wd0 \kern.2pt \box0}
\def\pmxbb#1{\setbox0\hbox{${#1}$} \copy0 \kern-\wd0 \kern.2pt
      \copy0 \kern-\wd0 \kern.2pt
      \copy0 \kern-\wd0 \kern.2pt \copy0 \kern-\wd0 \kern.2pt
      \copy0 \kern-\wd0 \kern.2pt \box0}

\def\cdotss{\mathinner{\cdotp\cdotp\cdotp\cdotp\cdotp\cdotp\cdotp
        \cdotp\cdotp\cdotp\cdotp\cdotp\cdotp\cdotp\cdotp\cdotp\cdotp
        \cdotp\cdotp\cdotp\cdotp\cdotp\cdotp\cdotp\cdotp\cdotp\cdotp
        \cdotp\cdotp\cdotp\cdotp\cdotp\cdotp\cdotp\cdotp\cdotp\cdotp}}

\font\frak=eufm10 scaled\magstep1 \font\fak=eufm10 scaled\magstep2 \font\fk=eufm10 scaled\magstep3
\font\scriptfrak=eufm10 \font\tenfrak=eufm10


\mathchardef\za="710B  
\mathchardef\zb="710C  
\mathchardef\zg="710D  
\mathchardef\zd="710E  
\mathchardef\zve="710F 
\mathchardef\zz="7110  
\mathchardef\zh="7111  
\mathchardef\zvy="7112 
\mathchardef\zi="7113  
\mathchardef\zk="7114  
\mathchardef\zl="7115  
\mathchardef\zm="7116  
\mathchardef\zn="7117  
\mathchardef\zx="7118  
\mathchardef\zp="7119  
\mathchardef\zr="711A  
\mathchardef\zs="711B  
\mathchardef\zt="711C  
\mathchardef\zu="711D  
\mathchardef\zvf="711E 
\mathchardef\zq="711F  
\mathchardef\zc="7120  
\mathchardef\zw="7121  
\mathchardef\ze="7122  
\mathchardef\zy="7123  
\mathchardef\zf="7124  
\mathchardef\zvr="7125 
\mathchardef\zvs="7126 
\mathchardef\zf="7127  
\mathchardef\zG="7000  
\mathchardef\zD="7001  
\mathchardef\zY="7002  
\mathchardef\zL="7003  
\mathchardef\zX="7004  
\mathchardef\zP="7005  
\mathchardef\zS="7006  
\mathchardef\zU="7007  
\mathchardef\zF="7008  
\mathchardef\zW="700A  

\newcommand{\be}{\begin{equation}}
\newcommand{\ee}{\end{equation}}
\newcommand{\ra}{\rightarrow}
\newcommand{\lra}{\longrightarrow}
\newcommand{\bea}{\begin{eqnarray}}
\newcommand{\eea}{\end{eqnarray}}
\newcommand{\beas}{\begin{eqnarray*}}
\newcommand{\eeas}{\end{eqnarray*}}
\def\*{{\textstyle *}}
\newcommand{\R}{{\mathbb R}}
\newcommand{\T}{{\mathbb T}}
\newcommand{\C}{{\mathbb C}}
\newcommand{\unit}{{\mathbf 1}}
\newcommand{\SL}{SL(2,\C)}
\newcommand{\Sl}{sl(2,\C)}
\newcommand{\SU}{SU(2)}
\newcommand{\su}{su(2)}
\def\ssT{\sss T}
\newcommand{\G}{{\goth g}}
\newcommand{\D}{{\rm d}}
\newcommand{\Df}{{\rm d}^\zF}
\newcommand{\de}{\,{\stackrel{\rm def}{=}}\,}
\newcommand{\we}{\wedge}
\newcommand{\nn}{\nonumber}
\newcommand{\ot}{\otimes}
\newcommand{\s}{{\textstyle *}}
\newcommand{\ts}{T^\s}
\newcommand{\oX}{\stackrel{o}{X}}
\newcommand{\oD}{\stackrel{o}{D}}
\newcommand{\obD}{\stackrel{o}{\bD}}
\newcommand{\pa}{\partial}
\newcommand{\ti}{\times}
\newcommand{\A}{{\cal A}}
\newcommand{\Li}{{\cal L}}
\newcommand{\ka}{\mathbb{K}}
\newcommand{\find}{\mid}
\newcommand{\ad}{{\rm ad}}
\newcommand{\rS}{]^{SN}}
\newcommand{\rb}{\}_P}
\newcommand{\p}{{\sf P}}
\newcommand{\h}{{\sf H}}
\newcommand{\X}{{\cal X}}
\newcommand{\I}{\,{\rm i}\,}
\newcommand{\rB}{]_P}
\newcommand{\Ll}{{\pounds}}
\def\lna{\lbrack\! \lbrack}
\def\rna{\rbrack\! \rbrack}
\def\rnaf{\rbrack\! \rbrack_\zF}
\def\rnah{\rbrack\! \rbrack\,\hat{}}
\def\lbo{{\lbrack\!\!\lbrack}}
\def\rbo{{\rbrack\!\!\rbrack}}
\def\lan{\langle}
\def\ran{\rangle}
\def\zT{{\cal T}}
\def\tU{\tilde U}
\def\ati{{\stackrel{a}{\times}}}
\def\sti{{\stackrel{sv}{\times}}}
\def\aot{{\stackrel{a}{\ot}}}
\def\sati{{\stackrel{sa}{\times}}}
\def\saop{{\stackrel{sa}{\op}}}
\def\bwa{{\stackrel{a}{\bigwedge}}}
\def\svop{{\stackrel{sv}{\oplus}}}
\def\saot{{\stackrel{sa}{\otimes}}}
\def\cti{{\stackrel{cv}{\times}}}
\def\cop{{\stackrel{cv}{\oplus}}}
\def\dra{{\stackrel{\xd}{\ra}}}
\def\bdra{{\stackrel{\bd}{\ra}}}
\def\bAff{\mathbf{Aff}}
\def\Aff{\sss{Aff}}
\def\bHom{\mathbf{Hom}}
\def\Hom{\sss{Hom}}
\def\bt{{\boxtimes}}
\def\sot{{\stackrel{sa}{\ot}}}
\def\bp{{\boxplus}}
\def\op{\oplus}
\def\bwak{{\stackrel{a}{\bigwedge}\!{}^k}}
\def\aop{{\stackrel{a}{\oplus}}}
\def\ix{\operatorname{i}}
\def\V{{\cal V}}
\def\cD{{\cal D}}
\def\cC{{\cal C}}
\def\cE{{\cal E}}
\def\cL{{\cal L}}
\def\cN{{\cal N}}
\def\cR{{\cal R}}
\def\cJ{{\cal J}}
\def\cT{{\cal T}}
\def\cH{{\cal H}}
\def\bA{\mathbf{A}}
\def\bI{\mathbf{I}}
\def\wh{\widehat}
\def\wt{\widetilde}
\def\ol{\overline}
\def\ul{\underline}
\def\Sec{\sss{Sec}}
\def\Lin{\sss{Lin}}
\def\ader{\sss{ADer}}
\def\ado{\sss{ADO^1}}
\def\adoo{\sss{ADO^0}}
\def\AS{\sss{AS}}
\def\bAS{\sss{AS}}
\def\bLS{\sss{LS}}
\def\bAP{\sss{AV}}
\def\bLP{\sss{LP}}
\def\AP{\sss{AP}}
\def\LP{\sss{LP}}
\def\LS{\sss{LS}}
\def\Z{\mathbf{Z}}
\def\oZ{\overline{\bZ}}
\def\oA{\overline{\bA}}
\def\cim{{C^\infty(M)}}
\def\de{{\cal D}^1}
\def\la{\langle}
\def\ran{\rangle}
\def\by{{\bi y}}
\def\bs{{\bi s}}
\def\bc{{\bi c}}
\def\bd{{\bi d}}
\def\bh{{\bi h}}
\def\bD{{\bi D}}
\def\bY{{\bi Y}}
\def\bX{{\bi X}}
\def\bL{{\bi L}}
\def\bV{{\bi V}}
\def\bW{{\bi W}}
\def\bS{{\bi S}}
\def\bT{{\bi T}}
\def\bC{{\bi C}}
\def\bE{{\bi E}}
\def\bF{{\bi F}}
\def\bP{{\bi P}}
\def\bp{{\bi p}}
\def\bz{{\bi z}}
\def\bZ{{\bi Z}}
\def\bq{{\bi q}}
\def\bQ{{\bi Q}}
\def\bx{{\bi x}}

\def\sA{{\sss A}}
\def\sC{{\sss C}}
\def\sD{{\sss D}}
\def\sG{{\sss G}}
\def\sH{{\sss H}}
\def\sI{{\sss I}}
\def\sJ{{\sss J}}
\def\sK{{\sss K}}
\def\sL{{\sss L}}
\def\sO{{\sss O}}
\def\sP{{\sss P}}
\def\sPh{{\sss P\sss h}}
\def\sT{{\sss T}}
\def\sV{{\sss V}}
\def\sR{{\sss R}}
\def\sS{{\sss S}}
\def\sE{{\sss E}}
\def\sF{{\sss F}}
\def\st{{\sss t}}
\def\sg{{\sss g}}
\def\sx{{\sss x}}
\def\sv{{\sss v}}
\def\sw{{\sss w}}
\def\sQ{{\sss Q}}
\def\sj{{\sss j}}
\def\sq{{\sss q}}
\def\xa{\tx{a}}
\def\xc{\tx{c}}
\def\xd{\tx{d}}
\def\xi{\tx{i}}
\def\xD{\tx{D}}
\def\xV{\tx{V}}
\def\xF{\tx{F}}
\def\dt{\xd_{\sss T}}
\def\vt{\textsf{v}_{\sss T}}
\def\vta{\operatorname{v}_\zt}
\def\vtb{\operatorname{v}_\zp}
\def\cM{\cal M}
\def\cN{\cal N}
\def\cD{\cal D}
\def\ug{\ul{\zg}}
\def\sTn{\stackrel{\scriptscriptstyle n}{\textstyle\sT}\!}
\def\sTd{\stackrel{\scriptscriptstyle 2}{\textstyle\sT}\!}
\def\stn{\stackrel{\scriptscriptstyle n}{\textstyle\st}\!}
\def\std{\stackrel{\scriptscriptstyle 2}{\textstyle\st}\!}


\setcounter{page}{1} \thispagestyle{empty}


\bigskip

\bigskip

\title{Lagrangian and Hamiltonian formalism in Field Theory:\\  a simple model
}

        \author{
        Katarzyna  Grabowska\thanks{The research financed by the Polish Ministry of Science and Higher Education under the
 grant N N201 365636.}\\
          {\it Physics Department}\\
                {\it University of Warsaw}}
\date{}
\maketitle
\begin{abstract}

The static of smooth maps from the two-dimensional disc to a smooth manifold can be regarded as a simplified
version of the Classical Field Theory. In this paper we construct the Tulczyjew triple for the problem and
describe the Lagrangian and Hamiltonian formalism. We outline also natural generalizations of this approach to
arbitrary dimensions.


\end{abstract}
\section{Introduction}The main purpose of this work is to implement the Tulczyjew  triple
approach of the Analytical Mechanics \cite{Tu1, Tu2} into the statics of multi-dimensional objects, i.e.
smooth maps from a disc $D\subset\R^n$ into a manifold $M$.  This problem can be regarded as a toy model for
the Classical Field Theory, since the set of smooth maps from $\R^n$ to $M$ can be treated as a set of
sections of the trivial bundle $pr_1:\R^n\times M\rightarrow \R^n$. In comparison with general geometric
approaches \cite{GIM1,GIM2,Tu0} the situation is considerably simplified, because the bundle is trivial and
the base manifold $\R^n$ has a canonical volume form and a canonical base of sections of the tangent bundle.
For $n=1$ and $M$ being the space of configurations of a mechanical system we recover the model of the
autonomous mechanics.

We work with this geometrically simple version of the Classical Field Theory to present the main ideas of our
approach to the Lagrangian and Hamiltonian formalism that differs from the ones which are present in the
literature \cite{K}. Since we skipped topological difficulties in this case, we could concentrate on the
recognition of physically important objects, like the phase space, phase dynamics, Legendre map, Hamiltonian,
etc. These issues are usually not elaborated well in the literature, as the Classical Field Theory models use
to concentrate on the Euler-Lagrange equations. Of course, we recover also the commonly accepted
Euler-Lagrange equations, this time without requiring any regularity of the Lagrangian.

The methods we use are based on expressing the theory in terms of differential relations rather than maps or
tensor fields. For the price of dealing with differential calculus of relations we get, in our opinion, better
understanding of geometric structures involved. It was also shown in \cite{GGU2, GG} that using the same
philosophy one can pass easily to the more complicated geometrical framework based on Lie or general
algebroids. In the case of Analytical Mechanics similar generalizations were proposed by many authors (e.g.
\cite{M1,M2,LMM}, but the approach presented in \cite{GGU2,GG}, being ideologically simpler, will be our
starting point.

We would like to point out that all the constructions we perform are motivated by the variational calculus
that we consider to be the fundamental idea of Classical Mechanics and Field Theory. The origin of geometric
structures we use lies in the rigorous formulation of the variational principle including boundary terms that
one can find in \cite{Tu0,NTu}. Nevertheless, we do not enter into details of the variational calculus and we
treat it rather as a guide-line for recognizing which geometrical structures are appropriate in this case.

The problem itself, i.e. the generalization of the symplectic framework for autonomous mechanics to higher
dimensions is not new and was first treated by G\"unther in \cite{Gun}. The underlying geometric structure of
G\"unther's theory, known as $k$-symplectic structure, was described systematically in \cite{Aw1, Aw2}.
Recently, Rey, Roman-Roy, Salgado and Valarino renewed the theory and described its Lie algebroid version
\cite{RRR}. Our work is also related to the multisymplectic approach to the Classical Field Theory developed
by Gotay, Isennberg, Marsden and others and presented e.g. in \cite{GIM1, GIM2, G1, G2}. The Tuczyjew triple
in the context of multisymplectic field theories appeared already in \cite{LMS}.

\medskip
For the presentation of our general idea, let us first recall the description of the dynamics of a classical
autonomous mechanical system without constraints. Let $M$ denote the manifold of positions of the system. The
trajectory is therefore a smooth path in $M$, i.e. a map from the time interval  $ [t_0, t_1]\subset \R$ into
$M$. We can try to describe our system in variational way, looking for those trajectories
$\gamma:\R\rightarrow M$ that, for the fixed time interval $[t_0, t_1]$, minimize the action functional
\be S(\gamma)=\int_{t_0}^{t_1}L(\st\gamma(t))\xd t.\ee
We assumed above that the Lagrangian is of first-order i.e. it is a function on the tangent bundle $\sT M$.
The curve $t\mapsto\st\gamma(t)$ will denote the {\it tangent prolongation} of the curve $\gamma$ in $M$. The
variational approach for the finite time interval leads to the Euler-Lagrange equations and the definition of
momenta. The space of momenta is usually called the {\it phase space} of the system. In the case of autonomous
mechanics, the phase space is just the cotangent bundle $\sT^\ast M$.

We describe the system by a first-order differential equation on the phase space, called the {\it phase
dynamics}. The phase dynamics $\mathcal{D}$ is described by a subset of $\sT\sT^\ast M$:
\be\mathcal{D}=\alpha_M^{-1}(\xd L(\sT M)),\ee
where $\alpha_M$ is the Tulczyjew isomorphism (defined in \cite{Tu1}) $\alpha_M: \sT\sT^\ast M\rightarrow
\sT^\ast\sT M$ and $\xd L(\sT M)$ is the image of the differential of the Lagrangian. A curve $t\mapsto
\zh(t)\in\sT^\ast M$ satisfies the phase dynamics if its tangent prolongation lies in $\mathcal{D}$. A curve
$t\mapsto\zg(t)\in\sT M$ satisfies, in turn, the corresponding Euler-Lagrange equation, if the curve
$t\mapsto\alpha_M^{-1}(\xd L(\zg(t)))\in\sT\sT^\ast M$ is the tangent prolongation of its projection to
$\sT^\ast M$ (see \cite{GGU2, GG}).

All the structures needed for generating the dynamics from the Lagrangian can be summarized in the following
commutative diagram:
\be\xymatrix@C-5pt{
& \sT\sT^\ast M \ar[rrr]^{\alpha_M} \ar[rdd]^{\sT\pi_M} \ar[ld]_{\zt_{\sT^\ast M}} & & & \sT^\ast\sT M
\ar[rdd]^{\pi_{\sT M}} \ar[ld]_{\zeta}\\
\sT^\ast M\ar[rrr]^/-25pt/{id}\ar[rdd]^/-20pt/{\pi_M} & & & \sT^\ast M\ar[rdd]^/-20pt/{\pi_M} &  \\
 & & \sT M\ar[rrr]^/-25pt/{id}\ar[ld]_{\zt_M} & & & \sT M\ar[ld]_{\zt_M} \\
&  M \ar[rrr]^{id}& & &  M }\ee The Lagrangian formulation can also be obtained from the variational
principle, when instead of finite domain of integration we use the so called de Rham current with one-point
support (see \cite{Tu0}).

It may happen that the phase dynamics is an implicit differential equation, i.e. it is not the image of a
vector field. In some cases, however, the phase dynamics is the image of a Hamiltonian vector field for some
function $H:\sT^\ast M\rightarrow \R$. In such a case we can write:
\be\mathcal{D}=\beta_M^{-1}(\xd H(\sT^\ast M)),\ee
where $\beta_M$ is the canonical isomorphism between $\sT\sT^\ast M$ and $\sT^\ast\sT^\ast M$ given by the
canonical symplectic form $\omega_M$ on $\sT^\ast M$,
$$\beta_M:\sT\sT^\ast M\longrightarrow\sT^\ast\sT^\ast M, \qquad\langle\beta_M(v),w\rangle=\omega_M(v,w).$$
Let us recall for the future reference that the canonical symplectic form $\omega_M$ is defined by
\begin{equation}\label{symp}\omega_M=\xd \vartheta_M,\end{equation}
where $\vartheta_M$ is the Liouville form given by
\begin{equation}\label{liouv}\vartheta_M(v)=\langle \zt_{\sT^\ast M}(v), \sT\pi_M(v) \rangle.\end{equation}

The structures needed for Hamiltonian mechanics can be presented in the following commutative diagram:
\be\xymatrix@C-5pt{
& \sT^\ast\sT^\ast M  \ar[rdd]^{\zx} \ar[ld]_{\pi_{\sT^\ast M}} & & & \sT\sT^\ast M\ar[lll]_{\beta_M}
\ar[rdd]^{\sT\pi_{M}} \ar[ld]_{\zt_{\sT^\ast M}}\\
\sT^\ast M\ar[rdd]^/-20pt/{\pi_M} & & & \sT^\ast M\ar[rdd]^/-20pt/{\pi_M}\ar[lll]_/-25pt/{id} &  \\
 & & \sT M\ar[ld]_{\zt_M} & & & \sT M\ar[ld]_{\zt_M}\ar[lll]_/-25pt/{id} \\
&  M & & &  M\ar[lll]_{id} }\ee

The formulation of the autonomous mechanics described above has at least two important features when compared
with the ones in textbooks: it is very simple and can be easily generalized to more complicated cases
including constraints, nonautonomous mechanics, and mechanics on algebroids \cite{GGU1, GGU2}. And last but
not least: we need no regularity conditions for the Lagrangian. The Lagrangian can be a function, but it can
be replaced by a family of functions generating a Lagrangian submanifold in $\sT^\ast\sT M$, as it happens in
the case of a relativistic particle in the Minkowski space \cite{TU}. The crucial role is played by two
mappings: $\alpha_M$ and $\beta_M$.


In what follows we replace `one dimensional' objects, like time intervals and paths in a manifold $M$ by `two
dimensional objects', like discs and maps $u:\R^2\rightarrow M$. We decided to keep $n=2$ just for simplicity.
However, generalization of our results to any natural $n$ is straightforward.

We shall then find the phase space and the analog of the fundamental map $\alpha_M$ that allows us to obtain
the phase equations from the Lagrangian. Then, we continue with the Hamiltonian formalism by recognizing what
kind of a geometric object the Hamiltonian is, and by finding an analog of the map $\beta_M$.

\section{Notation} Let $M$ be a smooth manifold of dimension $m$. We denote by $\zt_M \colon \sT M
\rightarrow M$ the tangent vector bundle and by $\zp_M:\sT^\ast M\rightarrow M$ the cotangent bundle of the
manifold $M$. If $(q^a)_{a=1}^m$ is a local coordinate system in $U\subset M$, then we have the induced
coordinate systems $(q^a, \dot q^b)$ in $\zt_M^{-1}(U)\subset\sT M$ and $(q^a, p_b)$ in
$\pi_M^{-1}(U)\subset\sT^\ast M$. The above coordinates correspond to local sections $(\xd q^b)$ of $\sT^\ast
M$ and $(\frac{\partial }{\partial q^b})$ of $\sT M$, respectively.

Let $u$ be a smooth map from $\R^2$ to $M$. Since in the source space $\R^2$ we have two distinguished vector
fields $\partial_{x^i}=\frac{\partial }{\partial x^i}$, $i=1,2$, the first jet $\sj^1 u(x)$ of the mapping $u$
at a point $x=(x^1, x^2)$ can be identified with a pair of vectors tangent to $M$ at the point $u(x)$, i.e.
$$\sj^1 u(x)=(v_1, v_2),\qquad\text{where}\qquad
v_i=(\sT_{x}u)(\partial_{x^i})\in\sT_{u(x)}M.
$$
Therefore the set $\sJ_{x}(\R^2, M)$ of the first jets of maps $\R^2\rightarrow M$ at $x=(x^1, x^2)\in\R^2$
will be identified with
$$\sTd M=\sT M\times_M\sT M$$
and the element of $\sTd M$, corresponding to $u$ at $x$, will be denoted by $\std u(x)$. The manifold $\sTd
M$ has a natural bundle structure over $M$:
$$ \sTd M\ni\;\std u(x)\longmapsto u(x)\in M.$$
The above projection will be denoted by $\zt^2_M$. Like for the tangent bundle, we have the adapted set of
local coordinates $(q^a, \dot q^{b}_1,\dot q^{c}_2)$ in $(\zt^2_M)^{-1}(U)\subset\sTd M$. The bundle $\zt^2_M$
is a vector bundle. Its dual $(\sTd M)^\ast$ can be identified with
$$\sTd\!{}^\ast M=\sT^\ast M\times_M\sT^\ast M$$
with the obvious pairing with $\sTd M$. The projection onto $M$ in the dual bundle will be denoted by
$\zp^2_M$. We have also the adapted set of local coordinates $(q^a,  q_{b\,1},p_{c\,2})$ in
$(\pi^2_M)^{-1}(U)\subset\sTd\!{}^\ast M$.

\section{Variational approach}
We start with Variational Calculus which is our guide-line for recognizing geometrical objects representing
physical quantities.  Let $L$ be a smooth function on the manifold $\sTd\! M$ of the first jets of maps from
$\mathcal{C}^\infty(\R^2,M)$; we will call $L$ a Lagrangian. Any Lagrangian defines an action functional $S$
on maps $u: D\rightarrow M$ from the unit disc $D\subset \R^2$ into M:
$$S(u)=\int_D L(u^a(x), \dot u^b_1(x), \dot u^c_2(x))\xd x^1\xd x^2,$$
where $u^a(x)=q^a(u(x))$ and $u^a_1(x)=\frac{\partial u^a}{\partial x^1}(x)=\dot q^a_1(\std u(x))$,
$u^a_2(x)=\frac{\partial u^a}{\partial x^2}(x)=\dot q^a_2(\std u(x))$. Note that the fact that Lagrangian can
be just a function on $\sTd M$ is due to the existence of the canonical volume form $\xd x^1\wedge \xd x^2$ on
$\R^2$. We can therefore identify scalar densities, i.e. objects that can be integrated, with functions.

Variations of $u$ are maps $\delta u$ from $D$ to $\sT M$ covering $u$:
$$\xymatrix{
D \ar[r]^{\delta u}\ar[dr]^u & \sT M\ar[d]^{\zt_M} \\
 & M
}$$ and coming from homotopies $\chi\in\mathcal{C}^\infty(\R\times\R^2,M)$ of maps
$\mathcal{C}^\infty(\R^2,M)$. If $u(x)=\chi(0, x)$, then $\delta u(x)$ is a vector tangent to the curve
$t\mapsto \chi(t, x)$ at $t=0$. In the following we perform the standard calculus of a variation of $S$ with
respect to the variation $\delta u$:
\beas
\langle \xd S, \delta u \rangle= \frac{\xd}{\xd t}_{|t=0}\int_D L(\chi^a(t,x), {\textstyle
\frac{\partial}{\partial x^1}} \chi^b(t,x),
{\textstyle \frac{\partial}{\partial x^2}} \chi^c(t,x))\xd x^1\xd x^2= \\
\int_D\left( \frac{\partial L}{\partial q^a}\delta u^a+ \frac{\partial L}{\partial \dot u^b_1}\delta \dot
u^b_1+ \frac{\partial L}{\partial \dot u^c_2}\delta \dot u^c_2 \right)\xd x^1\xd x^2=\ldots
\eeas
where $$\delta \dot u^a_i(x)=\frac{\partial  \chi^a}{\partial t\partial x_i}(0,x).$$ Using the Stokes theorem,
we obtain
$$\ldots=\int_D \left(
\frac{\partial L}{\partial q^a} -\frac{\partial}{\partial x_1}\frac{\partial L}{\partial q^a_1}
-\frac{\partial}{\partial x_2}\frac{\partial L}{\partial q^a_2}\right)\delta u^a\xd x^1\xd x^2 +\int_{\partial
D}\left( \frac{\partial L}{\partial q^a_1}\xd x_2-\frac{\partial L}{\partial q^a_2}\xd x_1\right)\delta u^a,
$$
where the last integral is calculated over $\partial D$ oriented as in the Stokes theorem, using the canonical
orien\-tation of $\R^2$. Looking for the stationary points of the action functional $S$ we put the condition
$\langle \xd S, \delta u \rangle=0$ for every $\delta u$, which means that
\be\frac{\partial L}{\partial q^a}-\frac{\partial}{\partial x_1}\left(\frac{\partial L}{\partial q^a_1}\right)
-\frac{\partial}{\partial x_2}\left(\frac{\partial L}{\partial q^a_2}\right)=0\label{el1}\ee on the disc and
$$\frac{\partial L}{\partial q^a_1}\xd x_2-\frac{\partial L}{\partial q^a_2}\xd x_1=0$$
on the boundary $\partial D$. The equation (\ref{el1}) is traditionally called the Euler-Lagrange equation.
The boundary term is an analog of the momentum in the Classical Mechanics. The momentum evaluated on a
variation $\delta u$ gives a one-form on $\R^2$ (to be integrated over $\partial D$). It follows that the
phase space is a space of covectors on $M$ with values in the cotangent space of $\R^2$. The cotangent bundle
of $\R^2$ is trivial, with the fiber being just $(\R^2)^\ast\simeq \R^2$, but we will keep the notation
$(\R^2)^\ast$ and use the basis $(\xd x^2, -\xd x^1)$ to identify it with $\R^2$. The phase space can be
therefore identified with
$$\sT^\ast M\otimes_M(\R^2)^\ast\simeq\sT^\ast M\otimes_M \R^2\simeq\sT^\ast M\times_M \sT^\ast M=\sTd^\ast M.$$
The Legendre map which associates a momentum to an infinitesimal configuration will be discussed later in the
section.

In the calculation of the differential of the action functional we have used implicitly a mapping
\be\kappa: \sT \!\sTd M\rightarrow \sTd \sT M\ee
defined as follows. Starting from a homotopy $\chi$ we can construct an element of $\sT\! \sTd M$ by taking
the tangent vector of the curve $\R\ni t\mapsto \std \chi(t,x)\in\sTd M$ at $t=0$. From the same homotopy  we
get $\R^2\ni x\mapsto \st\chi(t,x)\in\sT M$. The first jet at $x$ of the last map is an element of $\sTd \sT
M$. Therefore
\be\kappa(\st\std \chi(0,x))=\,\std\st\chi(0,x).\label{kappa1}\ee
The above definition is analogous to the definition of the canonical flip $\kappa_M:\sT\sT M\rightarrow\sT\sT
M$.

Using the local coordinate system $(q^a)$ on $M$, we can construct local coordinates on $\sTd\sT M$ and
$\sT\!\sTd M$. The coordinates on $\sTd\sT M$ will be denoted by
\be (q^a,\, \delta q^b,\, \dot q^c_1,\, \delta \dot q^d_1,\, \dot q^e_2,\, \delta \dot q^f_2)
\ee
and the ones on $\sT\!\sTd M$ by
\be (q^a,\, \dot q^b_1,\,  \dot q^c_2,\, \delta q^d,\,  \delta \dot q^e_1,\,  \delta \dot q^f_2).
\ee
Since
$$\kappa(q^a,\, \dot q^b_1,\,  \dot q^c_2,\, \delta q^d,\,  \delta \dot q^e_1,\,  \delta \dot q^f_2)=
(q^a,\, \delta q^b,\, \dot q^c_1,\, \delta \dot q^d_1,\, \dot q^e_2,\, \delta \dot q^f_2),$$ using the same
notation for coordinates in different spaces does not lead to any confusion.

\section{The Lagrangian side}

In the previous section we recognized the phase space as $\sTd\!{}^\ast M$. In the following it will be useful
to remember that the space $\sTd\!{}^\ast M$ has been identified with $\sT^\ast M\otimes_M(\R^2)^\ast$. W
shall define now a pairing between the space of jets of maps $p\in \mathcal{C}^\infty (\R^2, \sTd^\ast M)$ at
the point $(0,0)$, i.e. $\sTd\,\sTd\!{}^\ast M$, and jets of variations $\delta u\in\mathcal{C}^\infty(\R^2,
\sT M)$ at $(0,0)$, i.e. $\sTd \sT M$. (Using the structure of $\R^2$ we identify $\sJ^1(\R^2,\sTd\!{}^\ast
M)$ with $\R^2\times\sTd\,\sTd\!{}^\ast M$ and $\sJ^1(\R^2, \sT M)$ with $\R^2\times\sTd\sT M$, like in the
case of jets of maps $\mathcal{C}^\infty (\R^2, M)$).

The space $\sTd\,\sTd\!{}^\ast M$ has two vector bundle structures: the canonical one over $\sTd\!{}^\ast M$,
and the tangent one over $\sTd M$. The tangent projection $\sTd\pi^2_M$ can be constructed as follows: an
element $w\in\sTd\,\sTd\!{}^\ast M$ has a representative $\eta:\R^2\rightarrow\sTd\!^\ast M$, $w=\std
\eta(0,0)$. The projection $\sTd\pi^2_M$ is given by
$$\sTd\pi^2_M (w)=\std (\zp^2_M\circ \eta)(0,0).$$
Both bundles are vector bundles which form together a {\it double vector bundle} \cite{KU}:
\be\xymatrix@C-5pt{
 & \sTd\,\sTd\!{}^\ast M \ar[dl]_{\zt^2_{\sT^{2\ast} M}}\ar[rd]^{\sTd\pi^2_M} & \\
\sTd\!{}^\ast M\ar[rd]^{\pi^2_M} & & \sTd M \ar[dl]_{\zt^2_M} \\
 & M &
}\ee

Since the space of infinitesimal configurations $\sTd\! M$ has also a vector bundle structure over $M$, the
cotangent bundle $\sT^\ast\!\sTd M$ is a double vector bundle with the canonical projection on $\sTd M$ and
the second projection on $\sTd\!{}^\ast M$ (cf. for example \cite{GGU2,GU2}). The second projection is defined
as follows: any element $a$ of $\sT^\ast_{v}\!\sTd M$ is a linear function on the vector space $\sT_v\!\sTd
M$. Restricting $a$ to the subspace of vectors tangent to the fibre of projection $\zt^2_M$ (which is
isomorphic to the fibre itself), we can associate with it an element of $\sTd^\ast M$. The projection will be
denoted by $\zx$. The structure of the double vector bundle $\sT^\ast\!\sTd M$ can be put into the following
diagram:
\be\label{dvb1}\xymatrix@C-5pt{
 & \sT^\ast\sTd M \ar[dl]_{\zx}\ar[rd]^{\pi_{\sT^2 M}} & \\
\sTd\!{}^\ast M\ar[rd]^{\pi^2_M} & & \sTd M \ar[dl]_{\zt^2_M} \\
 & M &
}\ee The double vector bundle $\sT^\ast\!\sTd M$ is canonically isomorphic with $\sT^\ast\!\sTd\!{}^\ast M$.

Let $\std p$ and $\std \delta u$ be elements of $\sTd\,\sTd^\ast M$ and $\sTd\sT M$, respectively, such that
they have the same projection on $\sTd\! M$. Let $p$ and $\delta u$ denote the representatives covering the
same map $u:\R^2\rightarrow M$. Interpreting an element of $\sTd^\ast M$ as a covector on $M$ with values in
$(\R^2)^\ast$, we can define the mapping
$$\R^2\ni (x^1, x^2)\longmapsto \langle p(x^1, x^2), \delta u(x^1, x^2)\rangle \in (\R^2)^\ast\,,$$
where the target space is the fiber of $\sT^\ast\R^2$. The mapping can be viewed as a one-form on $\R^2$. The
differential of the above one-form is a two-form on $\R^2$ which, due to the existence of the canonical form
$\xd x^1\wedge\xd x^2$, can be identified with a function. The formula
\be\langle\!\langle \std p, \std \delta u \rangle\!\rangle\xd x^1\wedge \xd x^2=\xd \langle p, \delta u\rangle(0,0)\ee
defines a pairing  between $\sTd\,\sTd\!{}^\ast M$ and $\sTd\sT M$ over $\sTd\! M$. The above pairing is of
course degenerate.

Using the basis $(\xd x^2, -\xd x^1)$ of sections of $\sT^\ast\R^2$ and the local coordinate system $(q^a)$ on
$M$, we can construct local linear coordinates $(q^a, p^1_b, p^2_c)$ in the phase space $\sTd\!{}^\ast M$,
associated with local sections $(\xd q^b\otimes \xd x^2, -\xd q^c\otimes \xd x^1)$. Note that in our
convention the coordinates $(q^a, p^1_b, p^2_c)$ are associated with the the element $p^1_b\xd q^b\otimes \xd
x^2-p^2_c\xd q^c\otimes \xd x^1$. In the space of jets $\sTd\,\sTd\!{}^\ast M$ we have therefore the adapted
system of coordinates
$$(q^a, p^1_b, p^2_c,\; \dot q^d_1, \dot p^1_{e\,1}, \dot p^2_{f\,1},\; \dot q^g_2, \dot p^1_{h\,2}, \dot p^2_{k\,2} ).$$
On the other hand, in the space $\sTd\!\sT M$ of jets of variations we have also an adapted set of coordinates
induced by the coordinates on $M$:
$$(q^a, \delta q^b,\; \dot q^c_1, \delta \dot q^d_1,\; \dot q^e_2, \delta \dot q^f_2).$$
In the above coordinates the evaluation reads:
\be\langle\!\langle \std p, \std \delta u \rangle\!\rangle=(\dot p^1_{a\, 1}+\dot p^2_{a\, 2})\delta q^a+p^2_b\delta\dot q^b_2+
p^1_c\delta\dot q^c_1.\label{eval}\ee

Now we  are ready to define the mapping
$$\za:\sTd\,\sTd^\ast M\longrightarrow \sT^\ast\!\sTd M$$
by the condition
$$\langle \za(w), v \rangle=\langle\!\langle w, \kappa(v)\rangle\!\rangle, $$
where $w\in \sTd\,\sTd^\ast M$, $v\in\sT\!\sTd M$, and the evaluation on the left side of the equation is the
canonical evaluation between $\sT^\ast\sTd M$ and $\sT\sTd M$. The evaluation on the right side is the one
defined in (\ref{eval}). In local coordinates the mapping $\za$ reads
$$\za(q^a, p^1_b, p^2_c,\; \dot q^d_1, \dot p^1_{e\,1}, \dot p^2_{f\,1},\; \dot q^g_2, \dot p^1_{h\,2}, \dot p^2_{k\,2})=
(q^a, \dot q^d_1, \dot q^g_2,  \dot p^1_{f\,1}+ \dot p^2_{f\,2}, p^1_b, p^2_c).$$ The mapping $\za$ is an
analog of $\alpha_M:\sT\sT^\ast M\rightarrow\sT^\ast\sT M$ used by Tulczyjew in the autonomous mechanics.

\begin{definition} The {\it phase equations} for the system
with the Lagrangian $L$ are induced by the subset
\be \mathcal{D}=\za^{-1}(\xd L(\sTd M))\ee
in an obvious way: a mapping $p:\R^2\supset\mathcal{O}\rightarrow \sTd\!{}^\ast M$, where $\mathcal{O}$ is an
open subset of $\R^2$, is a solution of the phase equations if
$$\za(\std p(x,y))=\xd L(\;\std(\pi^2_M\circ p)(x,y)\;)$$
for any $(x,y)\in\mathcal{O}$.
\end{definition}

The important difference with the case of Classical Mechanics is that $\za$ is no longer an isomorphism,
therefore $\za^{-1}$ is a relation only, not a mapping. In local coordinates we obtain
\be \dot p^1_{a\,1}+ \dot p^2_{a\,2}=\frac{\partial L}{\partial q^a}, \quad p^1_b= \frac{\partial L}{\partial \dot q^b_1},
\quad p^2_c= \frac{\partial L}{\partial \dot q^c_2}.\ee The {\it Legendre map}, that associates a momentum to
an infinitesimal configuration, is defined as:
\be\lambda_L: \sTd M\longrightarrow \sTd^\ast M, \qquad \lambda_L=\xd L\circ \zx.\ee
In coordinates it reads
$$ \lambda_L(q^a, \dot q^b_1, \dot q^c_2)=\left(q^a, \frac{\partial L}{\partial \dot q^b_1}, \frac{\partial L}{\partial \dot q^c_2}\right)$$

The {\it Euler-Lagrange equations} for configurations
$$u:\R^2\ni(x^1, x^2)\longmapsto (q^a(x^1, x^2))\in M$$
can be formulated in the following way
\be \std \lambda_L=\za^{-1}\circ\xd L \ee
that in coordinates reads
\begin{align*}  & \displaystyle\dot q^a_1=\frac{\partial q^a}{\partial x^1}\,, \\
& \displaystyle \dot q^a_2=\frac{\partial q^a}{\partial x^2}\,, \\
& \displaystyle \frac{\partial L}{\partial q^a}- \frac{\partial}{\partial x^1} \frac{\partial L}{\partial\dot
q^a_1}- \frac{\partial}{\partial x^2} \frac{\partial L}{\partial\dot q^a_2}=0\,.
\end{align*}
The equations we obtained are in full agreement with equations commonly accepted in Classical Filed Theory the
theory (cf. \cite{GIM1,Kij}).

All the structure needed for generating the phase equations from the Lagrangian can be presented in the
following diagram:
\be\xymatrix@C-5pt{
& \sTd\,\sTd^\ast M \ar[rrr]^{\za} \ar[rdd]^/-20pt/{\sTd\pi^2_M} \ar[ld]_{\zt^2_{\sT^{2\ast} M}} & & &
\sT^\ast\sTd M
\ar[rdd]^{\pi_{\sT^2 M}} \ar[ld]_{\zx}\\
\sTd^\ast M\ar[rrr]^/-25pt/{id}\ar[rdd]^/-20pt/{\pi^2_M} & & & \sTd^\ast M\ar[rdd]^/-20pt/{\pi^2_M} &  \\
 & & \sTd M\ar[rrr]^/-25pt/{id}\ar[ld]_{\zt^2_M} & & & \sTd M\ar[ld]_{\zt^2_M} \\
&  M \ar[rrr]^{id}& & &  M }\ee

\section{The Hamiltonian side}

In the autonomous mechanics the basic structure is the canonical symplectic form $\omega_M$ on the cotangent
bundle being the phase space. Using the form $\omega_M$, we associate the Hamiltonian vector field to any
Hamiltonian -- a smooth function on the phase space. In our case, the phase space is not a symplectic manifold
any more, but still we can establish a correspondence between the cotangent bundle of the phase space and the
space $\sTd\,\sTd\!{}^\ast M$ of jets of maps from $\mathcal{C}^\infty(\R^2, \sTd\!{}^\ast M)$. We present
here  two ways of constructing the appropriate mapping.

The first method is based on the fact that the phase bundle is a vector bundle over $M$, so we have the
canonical antisymplectomorphism (cf. \cite{GU2,KU})
$$\sT^\ast\sTd^\ast M\simeq \sT^\ast \sTd M. $$
Denoting the above antisymplectomorphism by $\mathcal{R}$ and composing it with $\za$ we obtain
\be\label{beta1}\beta=\za\circ \mathcal{R} : \sTd\,\sTd^\ast M\longrightarrow \sT^\ast\sTd^\ast M\ee
Since $\mathcal{R}$ and $\za$ are double vector bundle morphisms, we obtain the following diagram for the
mapping $\beta$:
\be\xymatrix@C-5pt{
& \sT^\ast\sTd^\ast M \ar[rdd]^/-20pt/{\zz} \ar[ld]_{\zp_{\sT^{2\ast} M}} & & & \sTd\,\sTd^\ast M
\ar[rdd]^{\sTd\zp^2_M} \ar[ld]_{\zt^2_{\sT^{2\ast}M}}\ar[lll]_{\beta}\\
\sTd^\ast M\ar[rdd]^/-20pt/{\pi^2_M} & & & \sTd^\ast M\ar[lll]_/-25pt/{id}\ar[rdd]^/-20pt/{\pi^2_M} &  \\
 & & \sTd M\ar[ld]_{\zt^2_M} & & & \sTd M\ar[lll]_/-25pt/{id}\ar[ld]_{\zt^2_M} \\
&  M & & &  M\ar[lll]^{id} }\ee The projection $\zz: \sT^\ast\sTd^\ast M\rightarrow \sTd M$ is an analog of
$\zx$  in (\ref{dvb1}) and a particular instance of the canonical projection $\sT^\ast E^\ast\ra E$ for a
vector bundle $E$. If we use the system of local coordinates on $\sTd\,\sT^\ast M$ as previously,
$$(q^a, p^1_b, p^2_c,\; \dot q^d_1, \dot p^1_{e\,1}, \dot p^2_{f\,1},\; \dot q^g_2, \dot p^1_{h\,2}, \dot p^2_{k\,2} ),$$
and the system of coordinates
$$(q^a, p^1_b, p^2_c, \varphi_d, \psi_1^e, \psi_2^f ),$$
derived from the coordinates $(q^a, p^1_b, p^2_c)$ and associated to the local sections $(\xd q^a, \xd p^1_b,
\xd p^2_c)$, we get that
\be\label{beta}\beta(q^a, p^1_b, p^2_c,\; \dot q^d_1, \dot p^1_{e\,1}, \dot p^2_{f\,1},\; \dot q^g_2, \dot p^1_{h\,2}, \dot p^2_{k\,2} )=
(q^a, p^1_b, p^2_c, -\dot p^1_{d\,1}-\dot p^2_{d\,2}, \dot q^e_1, \dot q^f_2).\ee For any Hamiltonian
$$H:\;\sTd^\ast M\longrightarrow \R,$$
the phase dynamic is represented by the subset
\be\label{ham1} \mathcal{D}=\beta^{-1}(\xd H(\sTd\!{}^\ast M))\subset\sTd\,\sTd\!{}^\ast M\,.\ee
Also in this case, $\beta^{-1}$ is a relation only. In local coordinates we obtain the phase equations
\beas
-\dot p^1_{d\,1}-\dot p^2_{d\,2}=\frac{\partial H}{\partial q^d}, \\
\dot q^e_1=\frac{\partial H}{\partial p^1_e},\\
\dot q^f_2=\frac{\partial H}{\partial p^2_f}.\eeas

An alternative way of constructing the mapping $\beta$ does not refer to the map $\alpha$. Let us denote by
$pr_1$, $pr_2$ the projections on the first and the second factor of $\sTd^\ast M=\sT M\times_M\sT M$. In
local coordinates, we have
\beas pr_1: \sTd^\ast M\ni(q^a, p^1_b, p^2_c)\longmapsto (q^a, p^1_b)\in \sT^\ast M, \\
pr_2: \sTd^\ast M\ni(q^a, p^1_b, p^2_c)\longmapsto (q^a, p^2_b)\in \sT^\ast M.\eeas Applying the tangent lift
to the both projections we obtain
\beas \sT pr_1: \sT \sTd^\ast M\ni(q^a, p^1_b, p^2_c,\dot q^a, \dot p^1_b, \dot p^2_c)\longmapsto
(q^a, p^1_b, \dot q^a, \dot p^1_b )\in \sT\sT^\ast M, \\
\sT pr_2: \sT \sTd^\ast M\ni(q^a, p^1_b, p^2_c, \dot q^a, \dot p^1_b, \dot p^2_c)\longmapsto (q^a, p^2_b, \dot
q^a, \dot p^2_b)\in \sT \sT^\ast M.\eeas Composing the cartesian product of the above tangent mappings with
the inclusion
$$\imath: \sTd\,\sTd^\ast M\hookrightarrow\sT \sTd^\ast M\times \sT \sTd^\ast M,$$
we get
\begin{align*}(\sT pr_1\times \sT pr_2)\circ\imath:\sTd\,\sTd^\ast M&\longrightarrow \sT \sT^\ast M\times \sT \sT^\ast M\\
(q^a, p^1_b, p^2_c,\; \dot q^d_1, \dot p^1_{e\,1}, \dot p^2_{f\,1},\; \dot q^g_2, \dot p^1_{h\,2}, \dot
p^2_{k\,2} )\; &\longmapsto\; \left((q^a, p^1_b,\dot q^a_1, \dot p^1_{1b}),\; (q^a, p^2_b, \dot q^a_2, \dot
p^2_{2b})\right).\end{align*} To the both factors of the image of the composition $(\sT pr_1\times \sT
pr_2)\circ\imath$ we apply the canonical map $\beta_M:\sT\sT^\ast M\rightarrow \sT^\ast\sT^\ast M$ that comes
from the canonical symplectic form $\omega_M$ on the cotangent bundle $\sT^\ast M$. The target space of the
composition
$$(\beta_M\times\beta_M)\circ(\sT pr_1\times \sT pr_2)\circ\imath$$
is therefore
$$\sT^\ast\sT^\ast M\times \sT^\ast\sT^\ast M$$
which, in turn, can be mapped to $\sT^\ast \sTd^\ast M$ by means of the phase lift of the inclusion
$$\jmath:\sTd^\ast M\hookrightarrow\sT^\ast M\times\sT^\ast M.$$
Finally, we end up with the map
\be\label{beta2}\jmath^\ast\circ(\beta_M\times\beta_M)\circ(\sT pr_1\times \sT pr_2)\circ\imath:\;\;
\sTd\,\sTd^\ast M\longrightarrow \sT^\ast\sTd^\ast M.\ee

\begin{proposition} The mappings defined in (\ref{beta1}) and (\ref{beta2}) coincide, i.e.
$$\beta=\jmath^\ast\circ(\beta_M\times\beta_M)\circ(\sT pr_1\times \sT pr_1)\circ\imath.$$
\end{proposition}

\noindent{\bf Proof:} Let us start with recalling the definition of the canonical isomorphism $\mathcal{R}_E$
for a general vector bundle $E\rightarrow M$. The graph of $\mathcal{R}_E$ is the Lagrangian submanifold
generated in $\sT^\ast(E\times E^\ast)\simeq\sT^\ast E\times \sT^\ast E^\ast$ by the evaluation function
$$E\times_M E^\ast\ni (y,a)\longmapsto \langle a, y\rangle\in \R\,,$$
defined on the submanifold $E\times_M E^\ast$ of $E\times E^\ast$. We see that, by definition, for any element
$\varphi\in\sT^\ast E$, its image $\mathcal{R}_E(\varphi)$ has the same projections onto $E$ and $E^\ast$ as
$\varphi$. If we take now two curves
\begin{align*}
\R \ni t & \longmapsto \gamma(t)\in E, \\
\R \ni t & \longmapsto \eta(t)\in E^\ast,
\end{align*}
covering the same curve in $M$ and such that $\gamma(0)$ and $\eta(0)$ are equal to the projections of
$\varphi$ to $E$ and $E^\ast$ respectively, we can write
$$\langle (\varphi, \mathcal{R}_E(\varphi)), (\st\gamma(0), \st\eta(0))\rangle=
\frac{\xd}{\xd t}_{|t=0}\langle\eta(t), \gamma(t)\rangle,$$ or
\be \langle \mathcal{R}_E(\varphi), \st\eta(0)\rangle = \frac{\xd}{\xd t}_{|t=0}\langle\eta(t), \gamma(t)\rangle
-\langle \varphi, \st\gamma(0)\rangle.\label{er}\ee

Let now $\psi:\R^2\rightarrow \sT^\ast M$ be a homotopy such that $\psi(0,0)$ is the projection of $v$ and $w$
on $\sT^\ast M$, the curve $a\rightarrow \psi(a,0)$ is a representative of $v$, and $b\rightarrow \psi(0,b)$
is a representative of $w$. Using the definitions of $\omega_M$ and $\vartheta_M$ (see
(\ref{symp},\ref{liouv}), we can write that
\begin{multline*}\omega_M(v,w)=
\frac{\xd}{\xd b}_{|b=0}\theta_M(\st\psi(\cdot,b)(0))-
\frac{\xd}{\xd a}_{|a=0}\theta_M(\st\psi(a,\cdot)(0)=\\
\frac{\xd}{\xd b}_{|b=0}\langle \psi(0,b), \st(\pi_M\circ\psi(\cdot,b))(0)\rangle- \frac{\xd}{\xd
a}_{|a=0}\langle \psi(a,0), \st(\pi_M\circ\psi(a,\cdot))(0)\rangle.\end{multline*} We can simplify the above
formula a little  bit introducing curves
$$p_1(a)=\psi(a,0),\qquad p_2(b)=\psi(0,b),$$
which represent $v$ and $w$, respectively, and a homotopy in $M$ defined by
$$\rho=\pi_M\circ \psi.$$
In the new notation we have
\be\omega_M(v,w)=\frac{\xd}{\xd b}_{|t=0}\langle p_2(t), \st\rho(\cdot,b))(0)\rangle-
\frac{\xd}{\xd a}_{|a=0}\langle p_1(a), \st(\rho(a,\cdot))(0)\rangle.\label{omega}\ee

Now, we can start the main part of the proof, which is done by simple calculations. We are going to show that
the following diagram is commutative:
\be\xymatrix{ & \sTd\,\sTd^\ast M \ar[rd]^{\alpha}\ar[ld]_{\tilde\beta} & \\
\sT^\ast\!\sTd^\ast\! M & &\sT^\ast\!\sTd\! M \ar[ll]_{\mathcal{R}} }\label{proof1}\ee We have used the symbol
$\tilde\beta$ for the long expression $\jmath^\ast\circ(\beta_M\times\beta_M)\circ(\sT pr_1\times \sT
pr_1)\circ\imath$ and the shorter $\mathcal{R}$ for $\mathcal{R}_{\sT^2 M}$. The commutativity of the diagram
(\ref{proof1}) means that, for any $w\in \sTd\,\sTd^\ast M$ and $u\in \sT\sTd^\ast\! M$ such that they have
the same projection on $\sTd^\ast M$, we have
\be\langle\tilde\beta(w),u\rangle=\langle\mathcal{R}\circ\alpha(w),u\rangle.\label{proof2}\ee
To calculate the right-hand side of the equation (\ref{proof2}), we need the representatives of $u$ and $w$.
Let
$$\R\ni t\longmapsto \eta(t)=(\eta_1(t), \eta_2(t))\in \sTd^\ast M$$
be a curve that represents the vector $u$, i.e.
$$u=\st\eta(0).$$
The element $w$ of $\sTd\,\sTd^\ast M$ is represented by the mapping
$$\R^2\ni(x,y)\longmapsto p(x,y)=(p_1(x,y), p_2(x,y))\in \sTd^\ast M.$$
W can also choose a mapping
$$\R^3\ni(t,x,y)\longmapsto \chi(t,x,y)\in M$$
such that
$$\chi(t,0,0)=\pi^2_M\circ \eta(t),
\chi(0,x,y)=\pi^2_M\circ p(x,y).$$ Let us start with the R.H.S. of (\ref{proof2}) using (\ref{er}):
\begin{multline*}
\langle\,\mathcal{R}\circ\alpha(w),u\,\rangle=
\langle\,\mathcal{R}\circ\alpha(\,\std p(0,0)\,),\;\st\eta(0)\,\rangle= \\
\frac{\xd}{\xd t}_{|t=0}\langle\,\eta(t),\; \std\chi(t,\cdot, \cdot)(0,0)\,\rangle- \langle \,\alpha(\,\std
p(0,0)\,),\; \st\std\chi(\cdot,\cdot,\cdot)(0,0,0)\,\rangle.\end{multline*} Using the definition of $\alpha$
and the tangent evaluation, we can write the above expression as
$$\frac{\xd}{\xd t}_{|t=0}\langle\,\eta(t),\; \std\chi(t,\cdot, \cdot)(0,0)\,\rangle-
\langle\!\langle \,\std p(0,0),\; \std\st\chi(\cdot,\cdot,\cdot)(0,0,0)\,\rangle\!\rangle.$$ Finally,
replacing $\eta$ by the pair $(\eta_1,\eta_2)$ and using the definition of the evaluation
$\langle\!\langle\cdot,\cdot\rangle\!\rangle$, we obtain the final form
\begin{multline}\langle\,\mathcal{R}\circ\alpha(w),u\,\rangle=\frac{\xd}{\xd t}_{|t=0}\langle\,\eta_1(t),\; \st\chi(t,\cdot, 0)(0) \,\rangle+
   \frac{\xd}{\xd t}_{|t=0}\langle\,\eta_2(t),\; \st\chi(t,0, \cdot)(0) \,\rangle+ \\
   -\frac{\xd}{\xd x}_{|x=0}\langle\,p_1(x,0),\; \st\chi(\cdot,x,0)(0) \,\rangle-
   \frac{\xd}{\xd y}_{|y=0}\langle\,p_2(0,y),\; \st\chi(\cdot,0,y)(0) \,\rangle.\label{rhs}\end{multline}

Before we start with the L.H.S., let us note that (using representatives) we have
$$\imath(\std p(0,0))=(\st p_1(\cdot, 0)(0),\; \st p_2(\cdot, 0)(0),\;
\st p_1(0, \cdot)(0),\; \st p_2(0, \cdot)(0))$$ and
$$\sT pr_1\times\sT pr_2 (\st p_1(\cdot, 0)(0),\; \st p_2(\cdot, 0)(0),\;
\st p_1(0, \cdot)(0),\; \st p_2(0, \cdot)(0))= (\st p_1(\cdot, 0)(0),\;\st p_2(0, \cdot)(0)).$$ Now, we can do
our calculation as follows:
\begin{multline}\label{proof9}
\langle\, \jmath^\ast\circ(\beta_M\times\beta_M)\circ(\sT pr_1\times \sT pr_1)\circ\imath(w),\; u,\rangle=\\
\langle\, \jmath^\ast\circ(\beta_M\times\beta_M)\circ(\sT pr_1\times \sT pr_1)\circ\imath(\std p(0,0)),\; \st\eta(0),\rangle= \\
\langle\,(\beta_M\times\beta_M)(\st p_1(\cdot, 0)(0),\;\st p_2(0, \cdot)(0)),\;
\sT\jmath(\st\eta(0))\,\rangle.\end{multline} Taking into account that
$$ \sT\jmath(\st\eta(0))=(\st\eta_1(0), \st\eta_2(0)),$$
we can express (\ref{proof9}) as
\begin{multline*}
\langle\,\beta_M(\st p_1(\cdot, 0)(0)),\; \st\eta_1(0)\,\rangle+
\langle\,\beta_M(\st p_2(0, \cdot)(0)),\; \st\eta_2(0)\,\rangle= \\
\omega_M(\st p_1(\cdot, 0)(0),\; \st\eta_1(0))+\omega_M(\st p_2(0, \cdot)(0),\; \st\eta_2(0)).\end{multline*}
Now, we can use (\ref{omega}) for the final form of the L.H.S of equation (\ref{proof2}):
\begin{multline}\langle\tilde\beta(w), u\rangle
=\frac{\xd}{\xd b}_{|b=0}\langle\,\eta_1(b),\; \st\chi(b,\cdot,0)(0)\,\rangle-
\frac{\xd}{\xd a}_{|a=0}\langle\,p_1(a,0),\; \st\chi(\cdot,a,0)(0)\,\rangle \\
\frac{\xd}{\xd b}_{|b=0}\langle\,\eta_2(b),\; \st\chi(b,0,\cdot)(0)\,\rangle- \frac{\xd}{\xd
a}_{|a=0}\langle\,p_2(0,a),\; \st\chi(\cdot,0,a)(0)\,\rangle.\label{lhs}\end{multline} It is easy to see that
(\ref{rhs}) and (\ref{lhs}) coincide. $\Box$

The above proof is very technical, but tracing the calculations one can make at least one important
observation. In the final form of the R.H.S. and the L.H.S. one can see that momentum $(p_1, p_2)$ is always
evaluated on $\st\chi(\cdot,x,y)(0)$, i.e. on a tangent vector with respect to the first parameter $t$ of
$\chi$. The latter can be denoted by $\delta u$ and understood as a variation of the configuration $u$. In
(\ref{rhs}) it appears by definition from the tangent evaluation, but in (\ref{lhs}) it comes from the
calculation. The geometrical structure reflects therefore the idea that the momentum is to be evaluated on the
variation rather than on the infinitesimal configuration. In Classical Mechanics the difference is not
visible, because both, infinitesimal configurations and variations, are tangent vectors. In our case the
difference is visible but does not have much consequence, because the bundle of momenta is the dual vector
bundle of the bundle of infinitesimal configurations. In more general cases of Field Theory it is no longer
true. The Hamiltonian side is then more complicated.

In Classical Mechanics the mapping $\beta_M$ comes from the canonical symplectic form on the cotangent bundle
$\sT^\ast M$. We can therefore ask, whether the mapping $\beta$ we have just constructed is related to some
tensor field, which can be regarded as a canonical structure on the phase space. It is easy to see that,
indeed, the mapping $\beta$ is related to the field
$$\stackrel{\scriptscriptstyle 2}{\omega}_M\in\text{Sec}(\sTd\!{}^\ast\sTd\!{}^\ast M\otimes \sT^\ast\sTd\!{}^\ast M),$$
which in local coordinates reads
$$\stackrel{\scriptscriptstyle 2}{\omega}_M=\xd^1q^a\otimes\xd p^1{}_a+\xd^2 q^b\otimes \xd p^2{}_b
-\xd^1p^1{}_c\otimes \xd q^c-\xd^2p^2{}_d\otimes \xd q^d.$$ Here, $(\xd^i q^a, \xd^ip^1{}_b, \xd^ip^1{}_b )$
is the basis in the first (if $i=1$) or the second (if $i=2$) factor of $\sTd\!{}^\ast\sTd\!{}^\ast
M=\sT^\ast\sTd\!{}^\ast M\times_{\sT^2{}^\ast M}\sT^\ast\sTd\!{}^\ast M$. Since
$$\sTd\!{}^\ast\sTd\!{}^\ast M\otimes \sT^\ast\sTd\!{}^\ast M\simeq \sT^\ast\sTd\!{}^\ast M\otimes\sT^\ast\sTd\!{}^\ast M
\times_{\sT^2{}^\ast M} \sT^\ast\sTd\!{}^\ast M\otimes\sT^\ast\sTd\!{}^\ast M,$$ $\stackrel{\scriptscriptstyle
2}{\omega}_M$ can be viewed as a 'bi-form' or a pair of symplectic forms on the first and the second factor of
$\sT^\ast M$. This establishes a connection between the mapping $\beta$ and the poli-symplectic formalism of
G\"unther (\cite{Gun}).

\section{Example}

Since our model is a very simple and designed to study geometrical structures related to the Classical Field
Theory rather than to describe real physical systems, it is not easy to find physically important examples. In
the literature, one can find the so called {\it bosonic string theory}. There are two approaches to the
subject. In one of them, due to Polyakow \cite{Pol},  configurations are mappings from a two-dimensional
manifold $X$ of the string into the product of Minkowski space and the space of symmetric tensors on $X$. It
means that not only the space-time position of the string is subject to variations, but also the metric on the
string itself. In the simpler approach by Nambu \cite{Go, Na}, the metric on the string is fixed to be the
pull-back of the Minkowski metric by the string space-time configuration. In the Nambu approach we deal
therefore with mappings from the two-dimensional manifold to the Minkowski space. In our example we will use
the Nambu version with another simplification by taking $X=\R^2$.

The Minkowski space $(M,V, \eta)$ is a four-dimensional affine space with the model vector space $V$ equipped
with a bilinear symmetric form $\eta$ of signature $(+\,-\,-\,-)$. We will denote by $\tilde\eta$ the
associated self-adjoint map from $V$ to $V^\ast$. Using the affine structure, we can identify the tangent
bundle $\tau_M:\sT M\rightarrow M$ with the trivial bundle $M\times V\rightarrow M$, and the cotangent bundle
$\pi_M:\sT^\ast M\rightarrow M$ with the trivial bundle $M\times V^\ast\rightarrow M$. The spaces that appear
in the Lagrangian picture are therefore
\begin{align*}
& \sTd M=M\times V\times V\,, \\
& \sTd^\ast M=M\times V^\ast\times V^\ast\,, \\
& \sT^\ast\sTd M=(M\times V\times V)\times (V^\ast\times V^\ast\times V^\ast)\,, \\
& \sTd\,\sTd^\ast M=(M\times V^\ast\times V^\ast)\times (V\times V^\ast\times V^\ast)\times (V\times
V^\ast\times V^\ast)\,.
\end{align*}
The first jet of a mapping $u:\R^2\rightarrow M$ at the point $(x^1, x^2)$ is identified with a triple
$(q,v_1, v_2)$, where $q=u(x^1,x^2)$, $v_1$ is a vector tangent to the curve $t\mapsto u(x^1+t,x^2)$ at $t=0$,
and $v_2$ is a vector tangent to the curve $t\mapsto u(x^1,x^2+t)$ at $t=0$. The Lagrangian at the point
$\sj^1u$ is the scalar density associated to $u^\ast \eta$ which (after identification with the function on
$M\times V\times V$) gives
\begin{equation}L(q,v_1,v_2)=\sqrt{-\det g}\,,\end{equation}
where
$$g=\left[\begin{array}{cc}
\eta(v_1, v_1) & \eta(v_1, v_2) \\
\eta(v_1, v_2) & \eta(v_2, v_2)
\end{array}\right]\,.$$
The Lagrangian is defined on the open set of $M\times V\times V $, where the determinant of the matrix $g$ is
negative. Denoting with $(q, p^1, p^2)$ an element of $\sTd^\ast M=M\times V^\ast\times V^\ast$, and with $(q,
p^1, p^2, v_1, p^1{}_1, p^2{}_1,v_2, p^1{}_2, p^2{}_2)$ an element of $\sTd\,\sTd^\ast M$, we obtain the phase
equations
\begin{align}
& v_1=\frac{\xd q}{\xd x_1}, \quad v_2=\frac{\xd q}{\xd x_2}, \\
& p^1=\frac{1}{\sqrt{-\det g}}
\left[\eta(v_1, v_2)\tilde\eta(v_2)-\eta(v_2, v_2)\tilde\eta(v_1)\right], \\
& p^2=\frac{1}{\sqrt{-\det g}}
\left[\eta(v_1, v_2)\tilde\eta(v_1)-\eta(v_1, v_1)\tilde\eta(v_2)\right], \\
& p^1{}_1+p^2{}_2=0.
\end{align}
The Legendre map is in our example reversible, therefore we can express infinitesimal configurations in terms
of momenta:
\begin{align}
& v_1=\frac{-1}{\sqrt{-\det g}}
\left[\eta(p^1, p^2)\tilde\eta^{-1}(p^2)-\eta(p^2, p^2)\tilde\eta^{-1}(p^1)\right], \\
& v_2=\frac{-1}{\sqrt{-\det g}} \left[\eta(p^1, p^2)\tilde\eta^{-1}(p^1)-\eta(p^1,
p^1)\tilde\eta^{-1}(p^2)\right].
\end{align}
In the above formulae we used the same letter $\eta$ for the bilinear form associated to $\eta$ on the dual
side. The matrix $g$, in terms of momenta, takes the form
$$g=\left[\begin{array}{cc}
-\eta(p^2, p^2) & \eta(p^1, p^2) \\
\eta(p^1, p^2) & -\eta(p^1, p^1)
\end{array}\right].$$
Starting from the Hamiltonian
$$H(q,p^1,p^2)=-\sqrt{-\det g}\,,$$
we obtain the phase equations of the form
\begin{align}
& \frac{\xd q}{\xd x^1}=\frac{1}{\sqrt{-\det g}}
\left[\eta(p^1, p^2)\tilde\eta^{-1}(p^2)-\eta(p^2, p^2)\tilde\eta^{-1}(p^1)\right], \\
& \frac{\xd q}{\xd x^2}=\frac{1}{\sqrt{-\det g}}
\left[\eta(p^1, p^2)\tilde\eta^{-1}(p^1)-\eta(p^1, p^1)\tilde\eta^{-1}(p^2)\right], \\
& \frac{\xd p^1}{\xd x^1}+\frac{\xd p^2}{\xd x^2}=0\,,
\end{align}
which are of course the same as the phase equations generated by the Lagrangian description of the system. Let
us finish this example with writing down the fundamental maps $\alpha$ and $\beta$: On the Lagrangian side we
have
$$
\sT^\ast\sTd M=(M\times V\times V)\times (V^\ast\times V^\ast\times V^\ast)\,,$$ so that
$$\alpha(q, p^1, p^2, v_1, p^1{}_1, p^2{}_1,v_2, p^1{}_2, p^2{}_2)=(q,v_1, v_2, p^1{}_1+p^2{}_2, p^1, p^2 )\,.$$
On the Hamiltonian side, if we identify $\sT^\ast\sTd^\ast M$ with $(M\times V^\ast\times
V^\ast)\times(V^\ast\times V\times V)$, we obtain
$$\beta(q, p^1, p^2, v_1, p^1{}_1, p^2{}_1,v_2, p^1{}_2, p^2{}_2)=(q, p^1, p^2, -p^1{}_1-p^2{}_2, v_1, v_2).$$

\section{Conclusions}
We have presented a toy model of a Classical Field Theory to introduce main concepts of a new approach to
Lagrange and Hamilton formalisms. The starting point was the Tulczyjew triple in the Classical Mechanics,
generalized now to the case of fields. In this approach all main ingredients are present: starting with a
Lagrangian, not only the Euler-Lagrange field equation has been derived, but also the phase space and phase
dynamics have been recognized, together with the Legendre map and the Hamiltonian picture. The latter suggests
that momenta are dual rather to infinitesimal variations (displacements) than to infinitesimal configurations
(`velocities'). The main difference with respect to the classical situation is that, to construct the phase
dynamics, relations are used instead of mappings. This approach, presented here for maps from the disc into a
manifold, can be naturally generalized to sections of a fibration and to an `algebroid' setting as well. We
postpone these studies to a separate paper.

\end{document}